\documentclass{article}
\usepackage{amssymb,bbm}
\usepackage{amsfonts}
\usepackage{latexsym}
\usepackage{amsmath}
\usepackage{amsthm}
\usepackage{geometry}
\def\eb{\begin{eqnarray}}
\def\ee{\end{eqnarray}}

\def\non{\nonumber}
\def\wt{\widetilde}
\def\wh{\widehat}
\def\mc{\mathcal}
\def\mf{\mathfrak}
\def\bb{\mathbbm}

\def\g{\gamma}
\def\D{\mc{D}}
\def\tbf{\textbf}
\def\bp{\begin{proof}[{\bf Proof}]}

\def\deru{\frac{d}{du}}

\newtheorem{T1}{Theorem}

\begin{document}

\begin{titlepage}



\vspace{5 mm}

\begin{center}
{\Large\bf Spin chains from dynamical quadratic  algebras}

\vspace{10 mm}

{\bf Zolt\'an Nagy\footnote{e-mail: nagy@ptm.u-cergy.fr}, Jean Avan\footnote
{e-mail: avan@ptm.u-cergy.fr}}

\vspace{15mm}

{\it Laboratoire de Physique Th\'eorique et Mod\'elisation\\
     Universit\'e de Cergy-Pontoise (CNRS UMR 8089), 2 avenue Adolphe Chauvin,\\
 Pontoise
95302 Cergy-Pontoise cedex, France }
\vspace{5mm}

\end{center}

\vspace{18mm}

\begin{abstract}
\noindent We present a construction of integrable quantum spin chains where local spin-spin interactions are weighted by
``position''-dependent potential containing abelian non-local spin dependance.
This construction applies to  the previously defined three general
quadratic reflection-type algebras: respectively non-dynamical, semidynamical, fully dynamical.
\end{abstract}

\vfill

\end{titlepage}
\setcounter{footnote}{0}
\section{Introduction}
One major use of quantum groups and more generally  Hopf algebra type structures is the construction of exactly solvable
quantum spin chain type models.
The Hopf algebra structure enables one to build an $N$-site quantum transfer matrix; the existence of commuting traces then yields
commuting Hamiltonians inducing the relevant spin chain dynamics; finally the quantum group (or similar algebraic)
 structure allows for application
of algebraic or analytical methods (like Bethe Ansatz and Quantum Inverse Scattering) to obtain the eigenstates and observables of the
spin chain. For a recent account of these constructions and references see e.g. \cite{ACDFR}.

The particular choices of Hamiltonians yielding local spin interactions is also a crucial practical aspect of this procedure. Non-local
quantum spin chains may in principle be built but their actual resolution is generically much more difficult, 
although non-local quantum $N$-body
interacting systems of Ruijsenaars-Schneider spin type have been constructed \cite{ABB,KrZa}, resolution of which is connected to
the theory of special functions (Askey-Wilson polynomials, etc.).

Construction of local-interaction integrable spin chains was extended to the situation of reflexion (more generally soliton preserving
\cite{Skl}
 and soliton non-preserving \cite{DeVWo}) algebras, introduced by Cherednik \cite{Che}, Sklyanin, Kulish \cite{Skl,KuSkl,KuSa}
yielding open spin chain dynamics with integrable boundary conditions.

It is our purpose here to give an extensive description of the construction of locally-interacting (in a sense to be precised)
integrable spin chains derived from the comodule structure and quantum trace formulas for three generic quantum quadratic
algebras: the nondynamical \cite{KuSa,FrMa,DM,DKM,AD}; the semi-dynamical \cite{NAR,ArCheFr}; and fully dynamical
\cite{FHS,GZhZh,PPKu, AhnKoo,BP}  ``braiding''
algebras to be defined soon.

By generic we mean that we will not assume any particular relation  between the structure matrices of the quadratic structure, other than
whatever is induced by self-consistency of the algebra (Yang-Baxter, unitarity and locality requirements). Soliton-preserving
and soliton-non-preserving reflection algebras will only be particular cases of our derivation as we will comment later.

Specifically, we will construct homogeneous spin chain models, acting on spin quantum spaces of the form $W^{\otimes N}$ where
$W$ is an $n$-dimensional vector space and $W \otimes W$ is identified with the representation space of the structure matrices. 
It would be possible
to extend our construction to inhomogeneous situations, where one uses the fusion procedures on structure matrices \cite{KuSa,MeNe}
to obtain
different ``spin'' vector spaces $W^{(i)}$ at different spin sites $i \in \{1,\ldots,N\}$. Such a construction is given
e.g. in \cite{ACDFR}.

In the dynamical situation the eigenfunctions depend on a set of $n$ variables $x_i$ (same as $dim W$) on which the Hamiltonians
may or may not also act.
Locality properties of the spin chain Hamiltonians built here are characterized as follows: The Hamiltonians will be a sum over
spin sites (labeled by indices $i$ from $1$ to $N$ and respectively associated with the $i$-th copy of $W$ in the quantum space)
of operators written as direct products of two objects with distinct properties: one operator acting non-trivially on the sole
nearest-neighbor and next-to-nearest neighbor spin variables i.e. with site indices $j \in \{i-a_0, \ldots , i+a_0\}$ for $a_0 \leq 2$; and
one operator acting on the other spin variables , but belonging to an abelian enveloping algebra generated by  $N$ isomorphic
sets of $n$ commuting operators in an abelian algebra $\mf{h}$ associated with each spin site and such that the variables $x_i$ are 
identified as coordinates on $\mf{h}^{\ast}$. $W$ is here a diagonalizable $\mf{h}$-
module.
In the non-dynamical case these abelian operators reduce
to the identity and the resulting Hamiltonian is local in the exact sense. In the dynamical cases these operators are easily
computed and one can characterize the Hamiltonian as having a limited, abelian, non-local character..

We will now recall the definitions, notations and relevant properties of these three algebraic structures.


\section{Quadratic algebras}
\subsection{Non-dynamical exchange algebras}

The generators encapsulated in a matrix $T$ obey an exchange relation characterized by four matrices:
\eb
A_{12} T_{1} B_{12} T_{2} = T_{2}  C_{12}  T_{1}  D_{12}
\ee
with ``unitarity'' relations:
\eb
A_{12}=A_{21}^{-1}, \qquad B_{12}=C_{21}, \qquad D_{12}=D_{21}^{-1}
\ee
and cubic consistency relations of Yang-Baxter type describing the two ways of exchanging $T_1T_2T_3$ into $T_3 T_2 T_1$, assuming
associativity of the algebra.
Here as in all that follows, spectral parameter dependence is implicit in the labeling by indices, respectively $1$ and $2$.
We do not assume a priori that structure matrices depend only on specific combinations  (sum or difference) of the spectral parameters.
Stemming from the works in \cite{Che, Skl, KuSa} these algebraic strucutres were extended to general situations in \cite{FrMa}.
Their fusion or coproduct-like structure were investigated in \cite{KuSkl, KRS, MeNe, AD}; their universal interpretation as
originating from a Drinfel'd twist of a tensor product $H \otimes H$ of a Hopf algebra was recently proposed in \cite{DM, DKM} for
reflection algebras.

A dual exchange relation is naturally associated to the defining one \cite{KuSa,FrMa, AD}:
\eb
\left(A^{-1}_{12}\right)^{t_1t_2} \ K_1 \ \left(\left(B^{t_1}_{12}\right)^{-1}\right)^{t_2} \ K_2
 = K_2 \left(\left(C_{12}^{t_2}\right)^{-1}\right)^{t_1} \ K_1 \ \left(D_{12}^{t_1t_2}\right)^{-1}
\ee

A solution $K$ to the dual structure is indeed needed to
build naturally commuting traces which take the form \cite{Skl}:
\eb
\mathcal{H} =Tr K^{t} T
\ee
\subsection{Semi-dynamical exchange algebras}

The generators now obey an exchange algebra where the structure matrices $A,B,C,D$ and the generator matrix $T$ depend on a set of
parameters $\{\lambda\}$ interpreted as coordinates on the dual $\mf{h}^{\ast}$ of a subalgebra
$\mf{h}$ of the underlying Lie algebra. We will
consider only the case of an abelian $\mf{h}$ (non-abelian situations are considered e.g. in \cite{Ping}).
In addition we assume that the dimension of $\mf{h}$ and $\mf{h}^{\ast}$ is identical to the dimension of the vector space $W$
defining the representations of the structure matrices (as assumed in \cite{Et}).
As a consequence one can a) choose a basis of $W$ such that the generators
$\{h_i : i=1, \ldots,n\}$ are diagonal; and b) choose a basis of $\mf{h}$ such that the new basis vectors $h_i$ are identified with
the basis diagonal matrices $h_i=E_{ii}$ ($(E_{ii})_{jk}=\delta_{ij}\delta_{ik}$). The parameters $\{\lambda\}$ are correspondingly redefined.

The exchange algebra reads:
\eb \label{semidyn}
A_{12}(\lambda) T_{1}(\lambda) B_{12}(\lambda) T_{2}(\lambda+\g h_1) = T_{2}(\lambda)  C_{12}(\lambda)  T_{1}(\lambda+\g h_2)  D_{12}(\lambda)
\ee
where one formally denotes
\eb
f(\lambda+\gamma h)=  f(\{\lambda_i+\g h_i\}) = \sum_{m \geq 0}\frac{\g ^m}{m!} \sum_{i_1, \dots ,i_m = 1}^n
\frac {\partial^m f(\lambda)}{\partial \lambda_{i_1} \dots
\partial \lambda_{i_m}}  h_{i_1} \dots h_{i_m}
\ee
for any function $f$ on $\mathbb{C}^n$ with values in
$\mathsf{U}(\mathfrak{h})$.
The structure matrices $A,B,C,D$ obey once again Yang-Baxter type equations; in particular $A$ obeys a pure YB equation:
\eb
 A_{12}(\lambda) \ A_{13}(\lambda) \ A_{23}(\lambda) \ &=& \ A_{23}(\lambda) \ A_{13}(\lambda) \ A_{12}(\lambda)
\ee
while $D$ obeys the dynamical Gervais-Neveu-Felder equation \cite{GN,Fe}
\eb\label{eqD}
 D_{12}(\lambda+\gamma h_{3}) \ D_{13} \ D_{23}(\lambda+\gamma h_{1}) \ &=& \ D_{23} \ D_{13}(\lambda+\gamma h_2)
\ D_{12}
\ee
hence the denomination ``semi-dynamical''.

The first example of this algebra was given in \cite{ArCheFr,ArFr}. General study of this algebra, with complete elucidation of
fusion, comodule, dual algebra and quantum trace structure, is given in \cite{NAR,NADR}. Commuting traces are of the form:
\eb
\mathcal{H}=Tr K^t T e^{\g \D}
\ee
where $T$ obeys the exchange relation (\ref{semidyn}) and $K$ obeys the dual exchange relation (see \cite{nous}) and $\D$ is the 
shorthand notation for the difference operator $\sum h_i \partial_{\lambda_i}$.

\subsection{Fully dynamical exchange algebras}

In this case the exchange algebra reads :
\eb\label{ExRel}
A_{12}(\lambda) T_{1}(\lambda+\g h_2) B_{12}(\lambda) T_{2}(\lambda+\g h_1) = T_{2}(\lambda+\g h_1)  C_{12}(\lambda)  T_{1}(\lambda+\g h_2)  D_{12}(\lambda)
\ee
Structure matrices all obey Gervais-Neveu-Felder  equations (\ref{eqD}), hence the denomination ``fully dynamical''. The first example
of such an algebra is found in \cite{PPKu,BP,AhnKoo} and later in \cite{FHS,FHLS,GZhZh} as ``dynamical boundary algebra'',
albeit in the particular
``reflection'' situation where all four matrices are derived from
one single $R$-matrix. This situation was recently shown to originate from a similar Drinfel'd twist construction \cite{KuM},
applied to a more general algebraic structure (see e.g. \cite{Ping, Et}).



It is useful to develop somewhat more the general properties and formulas for this third case, so as to explicitely give the
generalized version of the construction presented in \cite{FHS,GZhZh}, particularly the commuting trace formula for a general algebra.

\subsubsection{Fusion and comodule structures}

Intuitively, both fusion and comodule structures are related to tensor product. Whereas fusion operates on auxiliary spaces,
comodule operates on
quantum spaces, thus one obtains the extension of Sklyanin's double-row transfer matrix to the dynamical case. The fusion
structures are described in \cite{nous}. The comodule structure is characterized by:

\begin{T1}
Let $T_{q}$ be a representation of the algebra defined by (\ref{ExRel}).
Let $L_{q'}, R_{q'}$ denote a representation on another Hilbert space $\mc{H}_{q'}$ of the following exchange relations.
\eb \non
A_{12}(\lambda) L_{1q'}(u_1;\lambda+\g h_2) L_{2q'}(u_2;\lambda) & =& L_{2q'}(u_2;\lambda+\g h_1) L_{1q'}(u_1;\lambda) A_{12}(\lambda+\g h_{q'}) \\ \non
R_{1q'}(u_1;\lambda+\g h_2) B_{12}(\lambda) L_{2q'}(u_2;\lambda+\g h_1) & = & L_{2q'}(u_2;\lambda) B_{12}(\lambda+\g h_{q'}) R_{1q'}(u_1;\lambda) \\ \non
L_{1q}(u_1;\lambda) C_{12}(\lambda+\g h_{q'}) R_{2q'}(u_2;\lambda) & = & R_{2q'}(u_2;\lambda+\g h_2) C_{12}(\lambda) L_{1q'}(u_1;\lambda+\g h_2) \\ \non
D_{12}(\lambda+\g h_{q'}) R_{1q'}(u_1;\lambda) R_{2q'}(u_2;\lambda+\g h_1) & = & R_{2q'}(u_2;\lambda) R_{1q'}(u_1;\lambda+\g h_2) D_{12}(\lambda)
\ee
then
\eb
T_{1,qq'}(u_1;\lambda)= L_{1q'}(u_1;\lambda) T_{1q}(u_1;\lambda+\g h_{q'}) R_{1q'}(u_1;\lambda)
\ee
is a representation on $\mc{H}_q \otimes \mc{H}_{q'}$ of the original exchange algebra (\ref{ExRel}).
\end{T1}
\bp Straightforward by plugging the new $T_{1,qq'}$ into the original exchange relation and using the defining relations of $L$ and $R$.
\end{proof}

A simple example to this comodule algebra is provided on $\mc{H}_{q'}=W$ by the structure matrices. Namely,
$L_{1q'}(u_1)=A_{12}(u_1,\alpha;\lambda), R_{1q'}(u_1;\lambda)=B_{12}(u_1,\alpha;\lambda)$ or
$L_{1q'}(u_1)=C_{12}(u_1,\alpha;\lambda), R_{1q'}(u_1;\lambda)=D_{12}(u_1,\alpha;\lambda)$ for any $\alpha \in \mathbb{C}$.




Commuting traces can be constructed for this algebra as well. They take the form:
\eb\label{TMfuldyn}
\mc{H}=Tr e^{-\g\D} T e^{\g\D} (K^{SC}) ^{t}
\ee
where $(M^{SC})_{ij}=M_{ij}(\{x_k-\g \delta_{jk}\})$).For the proof of commutation and the definition of 
the corresponding dual algebra see \cite{nous,YaSa}.

Let us comment here on how this transfer matrix formula is related to the one in \cite{FHLS} initially written for the
$A_{n-1}^{(1)}$ IRF model.
The transfer matrix (\ref{TMfuldyn}) can be detailed as:
\eb
H(u,\lambda)=\sum_{ij} e^{-\partial_i} T_{ij} e^{\partial_j} e^{-\partial_j}K_{ij} e^{\partial_j}= \sum T^{-SL}_{ij} K^{-SL}_{ij}
e^{-\partial_i+\partial_j}
\ee
This operator acts on the Hilbert space of $V$-valued functions, i.e. on sums of $f(\lambda) \otimes |spin\rangle$ type tensor products.
It is easy to see that if $\sum T^{-SL}_{ij} K^{-SL}_{ij}$ is \textit{independent} of $\lambda$ then $H$ leaves the subspace generated by
$constant \otimes |spin\rangle$ invariant and can be restricted to it. The restriction of $H$ to this subspace is equal to
the restriction of $\tilde{H}=\sum T^{-SL}_{ij} K^{-SL}_{ij}$ and commutes at different values of the spectral parameter.
 Now, using the independence of $\tilde{H}$ on $\lambda$ we can lift its action back to the whole Hilbert
space again (which is a direct sum $\bigoplus_n \mathbb{C}\lambda^n \otimes |spin\rangle$)
and get the result $\left[\tilde{H}(u),\tilde{H}(v)\right]=0$.

The IRF case described in \cite{FHLS} falls precisely in this category as its transfer matrix is a rewriting
(via face-vertex intertwiners \cite{JiMiOk}) of a vertex
type transfer matrix and is thus independent of the face parameter $\lambda$.

\section{Spin chains}

The construction of local spin chains starts with the building of a transfer matrix; this procedure here exploits the comodule structure
in the following way: one starts with a scalar solution of the exchange relation and adds sites to build a  double-row transfer
matrix accordingly to the comodule structure. The transfer matrix is then closed on the other end of the chain by a scalar
representation of the ``dual'' exchange relation. We will not consider situations where one uses quantum representations
(also called dynamical) as in e.g. \cite{BaDo})
In all quadratic cases at least one particular representation of the comodule is given by the structure matrices and we will restrict our
analysis to this case.
One then takes the (logarithmic) derivative of the transfer matrix with respect to the auxiliary spectral parameter
at a particular value of the auxiliary and quantum spectral parameters, once
again yielding commuting Hamiltonians.
Locality is assured when certain structure matrices become a \emph{permutation} at the value $0$ of
the spectral parameters.
Indeed this allows for a decoupling inside the comodule structure, and a simplification of the
Hamiltonians down to terms involving only neighboring structure matrices.

Let us immediately comment on the essential difference between our construction of spin chains and the ones using similar algebraic
structures in \cite{ABB} or \cite{GZhZh}:

In \cite{ABB} dynamical quantum groups and the corresponding commuting traces were used to build Calogero and
Ruijsenaars type Hamiltonians. The procedure used the $u \to \infty$ limit of $R$ matrices $R_{12}(u=u_1-u_2)$ in order
to build a nonlocal spin interaction model where ``position'' and spin variables were tied together by the further required
assumption $N$(site number) $= n$ (matrix dimension).
A further reduction can eliminate the spin variables altogether  and lead to a scalar Ruijsenaars-Schneider
Hamiltonian.
In \cite{GZhZh} the Gaudin-type Hamiltonians are built by identifying the specific fixed values of the quantum space spectral
parameters (in a comodule construction of transfer matrix) with arbitrary ``positions'' of the spin sites. Locality of the
Hamiltonian is ensured by a semi-classical limit procedure.
By contrast, in this paper, we use matrices whose
limit $u_1,u_2\to 0$ yields a permutation, thus immediately giving rise to local Hamiltonians (in the sense characterized in the introduction).

Depending on the behavior at spectral parameters $u_1,u_2=0$ of
the structure matrices $A,B,C,D$, two cases can be distinguished.

\begin{itemize}
\item
Soliton-preserving case: Every structure matrix equals $P_{12}$ at spectral parameters $u_1,u_2=0$.
\item
Soliton-non-preserving case: $A_{12}(0,0)=D_{12}(0,0)=P_{12}$. $B$ and $C$ not necessarily have the same limit. In this case an alternating
spin chain should in general be considered in order to ensure locality.
However, at least in the case where $B(0,0)= \beta_2 P_{12} \delta_2$ it is not necessary to
define an alternating spin chain transfer matrix. The spin chain resulting from a non-alternating transfer matrix will not be local but
will be related to a local one by a product of $\beta$'s and $\delta$'s. This is the case in particular after redefinition of
the structure matrices (cf. below).

\end{itemize}

Let us now suppose that there exists a scalar representation $T$ and $K$ for the original and the dual exchange relation. This is in fact a
necessary additional input because the structure matrices themselves do not automatically provide a representation - unlike in the
$RTT$ case. Once these representations are given one can assure by redefinition of the structure matrices that one of them, say, $K$
is trivial. This is achieved by introducing the notion of isomorphic $K$-conjugate dual quadratic algebra:

We consider first the nondynamical case \cite{FrMa}. Let  $\chi$
be such a representation of the nondynamical dual algebra, that is
\eb \non (A_{12}^{-1})^{t_{12}} \chi_1 ((B_{12}^{t_1})^{-1})^{t_2}
\chi_2 = \chi_2 ((C_{12}^{t_2})^{-1})^{t_1} \chi_1
(D_{12}^{-1})^{t_{12}} \ee Using this scalar solution one defines
the $\chi$-conjugated dual algebra with the following structure
matrices and with the redefined $K$-matrix: \eb \non
&\tilde{A}^{dual}_{12}=\chi_1^{-1}\chi_2 ^{-1} A_{12}^{dual} \chi_1\chi_2, \quad \tilde{B}^{dual}_{12}=\chi_2^{-1} B^{dual}_{12} \chi_2& \\
&\tilde{C}^{dual}_{12}=\chi_1^{-1} C^{dual}_{12} \chi_1, \quad
\tilde{D}^{dual}_{12}=D_{12}^{dual} \quad \tilde{K}=\chi^{-1}K&
\non \ee

This $\chi$-conjugated dual  algebra is now endowed with a trivial
representation $\tilde{K}=\mathbbm{1} \otimes \mathbbm{1}$. Of
course, the redefinition of the dual matrices entails in turn the
redefinition of the structure matrices of the algebra itself as:
\eb & \tilde{A}_{12}=\chi_1^{t}\chi_2 ^{t} A_{12} (\chi_1^{-1})^{
t}(\chi_2^{-1})^{ t} \quad  \tilde{B}_{12}=\chi_2^t B_{12}
(\chi_2^t)^{-1}& \\
& \tilde{C}_{12}=\chi_1^{t} C_{12} (\chi_1^t)^{-1} \quad \tilde{D}_{12}=D_{12}& \non \ee obeying Yang-Baxter
equations. The new exchange relations are obeyed by the redefined $\tilde{T}=\chi^t T$ matrix.

Note also that unless a ``crossing'' relation of the type $(B_{12}^{t_2})^{-1}=(B_{12}^{-1})^{t_2}$, and similarly
for $C$, is satisfied by $B$ and $C$, the existence of a solution $\chi$ for the dual algebra does not imply the existence
of a solution related to $\chi$ for the original algebra.

For the dual semidynamical algebra the redefinitions are:
\eb
&\tilde{A}^{dual}_{12}=\chi_1^{-1}\chi_2 ^{-1} A_{12}^{dual} \chi_1\chi_2 \quad
\tilde{B}^{dual}_{12}=\chi_2^{-1} B^{dual}_{12} \chi_2(h_1) & \non \\
& \tilde{C}^{dual}_{12}=\chi_1^{-1} C^{dual}_{12} \chi_1(h_2) \quad \tilde{D}^{dual}_{12}=D_{12}^{dual}\quad
\tilde{K}=\chi^{-1}K & \non \ee
And for the algebra itself, the structure matrices become: \eb
& \wt{A}_{12}=\chi_1^t\chi_2 ^t A_{12} (\chi_1^{-1})^t (\chi_2^{-1})^t \quad  \wt{B}_{12}=\chi_2^t B_{12} (\chi_2^{-1})^t(h_1) \non & \\
& \wt{C}_{12}=\chi_1^t C_{12} (\chi_1^{-1})^t(h_2) \quad \wt{D}_{12}=D_{12} \non \quad \tilde{T}=\chi^t T&
\ee Note that in this case the above mentioned crossing relation is satisfied because of the partial zero
weight condition on $B$ and $C$ (discussed in \cite{NAR}) and that these partial zero weight properties are conserved by the
redefinition.

In the dynamical case the redefinitions for the dual algebra are:
\eb
& \tilde{A}^{dual}_{12}=\chi_1^{-1}\chi_2^{-1}(h_1) A_{12}^{dual} \chi_1(h_2)\chi_2 \quad
\tilde{B}^{dual}_{12}=\chi_2^{-1} B^{dual}_{12} \chi_2(h_1)& \non \\
& \tilde{C}^{dual}_{12}=\chi_1^{-1} C^{dual}_{12} \chi_1(h_2) \quad \tilde{D}^{dual}_{12}=D_{12}^{dual} \quad
\tilde{K}=\chi^{-1}K &\non \ee

And for the algebra itself:
\eb
& \wt{A}_{12}=\chi_1^t\chi_2 ^t(h_1) A_{12} (\chi_1^{-1})^t (\chi_2^{-1})^t(h_1) \quad  \wt{B}_{12}=\chi_2^t B_{12} (\chi_2^{-1})^t(h_1)
\non & \\
& \wt{C}_{12}=\chi_1^t C_{12} (\chi_1^{-1})^t(h_2) \quad \wt{D}_{12}=D_{12} \quad \tilde{T}=\chi^t T\non &
\ee Once again the existence of a solution $\chi$ for the dual algebra structure does not imply in general
the existence of a solution for the direct algebra, as commented in \cite{YaSa}. However, in a case where a
face-vertex correspondance can be established such a duality exists \cite{YaSa}.

\noindent\textbf{Important remarks:}  One has to assume  in the dynamical  case here that the particular solution $\chi$ is of zero weight,
i. e. diagonal
in our choice of basis for $\mf{h}$. The general case ($\chi$ non-diagonal) must be dealt with separately.

In all cases, this redefinition may, however, change the nature of the matrices $B,C$ from soliton preserving to soliton
non preserving, but one will
be in the situation evoked in the beginning as $B \sim \g_2 P_{12} \delta_2$.

We will see that in all cases the transfer matrix, and hence the Hamiltonians, obtained from the $\chi$-conjugate algebraic formulation
with $\chi=1$, will be exact conjugates (by some product of $\chi$'s) of the original transfer matrix with $\chi \neq 1$.  Hence,
by considering both SP and SNP situations in all three cases with $\chi =1$, plus the extra non-diagonal fully dynamical case,
one will indeed cover all possibilities for the form of the
Hamiltonians, up to conjugations.


\subsection{Nondynamical spin chain}

The \emph{soliton preserving} case is a simple extension of Skylanin's original construction of integrable open spin chains \cite{Skl}.
The transfer matrix takes
the form:
\eb
&&t(u)=tr_0 \chi_0^t(u) A_{0N}(u,u_N)\ldots A_{01}(u,u_1) T_0(u) \non \\
&&B_{01}(u,u_1) \ldots B_{0N}(u,u_N)
\ee
where $\chi$ is a scalar ($c$-number) solution of the dual exchange algebra.
This transfer matrix is related by conjugation to the  transfer matrix $\tilde{t}(u)$ with $\chi=1 $
obtained from the \emph{conjugated} algebra as follows:
\eb
&&\tilde{t}(u)=tr_0 \tilde{A}_{0N}(u,u_N)\ldots \tilde{A}_{01}(u,u_1) \tilde{T}_0(u) \non \\
&&\tilde{B}_{01}(u,u_1) \ldots \tilde{B}_{0N}(u,u_N)=\prod_i \chi_i^t \ t(u) \ \prod_i (\chi_i^{t})^{-1}  \non
\ee

Let us now take $\chi_0=1$. The Hamiltonian has the form:
\eb
\frac{d}{du}\bigg|_{\substack{u=0\\ u_i=0 }}\log t(u)&=& \sum_{j=1}^{N-1} \check{A}'_{j,j+1} +  \sum_{j=2}^{N-1} \check{B}'_{j+1,j}  +
T_1 \check{B}'_{21} T_1^{-1} +\non \\
&&tr_0 \check{A}'_{N0} + tr_0 \check{B}'_{0N} +T'_1T^{-1}_1 \non
\ee
where $\check{A}_{ij}=P_{ij}A_{ij}$ and all spectral parameters are set to $0$. Unlike Sklyanin in \cite{Skl}, we take the logarithmic
derivative here, because we do not assume $T(0)=1$.

The \emph{soliton nonpreserving} case is a generalization of the results in \cite{An}.

The transfer matrix resulting from the comodule reads:
\eb
&&t(u)=tr_0 \ \chi_0^t A_{02N}(u,u_{2N})C_{02N-1}(u,u_{2N-1}) \ldots A_{02}(u,u_2)C_{01}(u,u_1)T_0(u) \non\\
&&D_{01}(u,u_1)B_{02}(u,u_2)\ldots D_{02N-1}(u,u_{2N-1})B_{02N}(u,u_{2N})
\ee
This transfer matrix is related by a conjugation to the one obtained from the conjugated algebra as follows:
\eb
\tilde{t}(u)= \prod_{i \textit{ even}} \chi_i^t \ t(u) \ \prod_{i \textit{ even}} (\chi_i^t)^{-1}
\ee
Once again we continue working with $\chi=1$.
At $u=u_i=0$ this gives:
\eb
&&t(0)= C_{2N,2N-1}C_{2N-2,2N-3}\ldots C_{21} P_{24}P_{46}\ldots P_{2N-2,2N}P_{2N-1,2N-3}
\ldots P_{53}P_{31} \non \\
&&B_{32} \ldots B_{2N-1,2N-2} T_{2N}P_{12N}tr_0 P_{02N}B_{02N} \non
\ee
Let us define $X_k=tr_0 P_{0k}B_{0k}$ which acts on the $k$th site of the chain. Let us also suppose that $X$ is invertible: this
property is needed when calculating the Hamiltonian.

 All spectral parameters are set to $0$, because the regularity condition on $A$ fixes the $u_{2k}$'s and
the regularity condition on $D$ fixes the $u_{2k+1}$'s. So from now on spectral parameters will be omitted.
The logarithmic derivative with respect to the first spectral parameter of an operator in $End(V)\otimes End(V)$ depending
on two spectral parameters will be denoted by the corresponding curly letter, for example:
\eb
\mc{C}_{ij}:= \frac{d}{du}\bigg|_{u=0}C_{ij}(u,0) C_{ij}^{-1}(0,0)
\ee
and the conjugation will be denoted by:
\eb
Ad(A)\cdot B:=ABA^{-1}
\ee
The calculations yield the following Hamiltonian:
\eb
&&\frac{d}{du}\bigg|_{\substack{u=0\\ u_i=0 }}\log t(u)= \sum_{j=1}^N \mc{C}_{2j,2j-1}+  \sum_{j=1}^{N-1} Ad(C_{2j+2,2j+1})\cdot
\mc{\check{A}}_{2j,2j+2} +\non \\
&& \sum_{j=1}^{N-1} Ad(C_{2j+2,2j+1} C_{2j,2j-1})\cdot \mc{B}_{2j-1,2j+2} +\non \\
&& \sum_{j=1}^{N-1} Ad(C_{2j+2,2j+1}C_{2j,2j-1} B_{2j-1,2j+2})\cdot \mc{\check{D}}_{2j+1,2j-1}+
\non \\
&& tr_0 \ \check{A}'_{2N0} C_{2N,2N-1} P_{0,2N-1} B_{0,2N-1} X^{-1}_{2N-1} C^{-1}_{2N,2N-1} + \non \\
&& C_{2N,2N-1}\deru \log  X_{2N-1}(u)C^{-1}_{2N,2N-1} +  Ad(C_{21}T_2)\cdot \mc{\check{D}}_{12}+ C_{21} T'_2 T_2^{-1} C_{21}^{-1}\non
\ee

\subsection{Semidynamical spin chain}

Note that the weight conditions on $B$ and $C$ do not allow for soliton preserving boundary conditions.

Here, as always, the dynamical parameters will be omitted, only the shifts will be explicitely noted. $A_{i_1,\ldots, i_m}(h_{<}^{odd})$
stands for $A_{i_1, \ldots,i_m}$ shifted on all the spaces with odd indices greater than the indices of $A$, if any.

We start here to build the chain as $\ldots C_{01}T_0 D_{01} \ldots $. 
Notice that the reason for a shift on the sole odd spaces lies indeed in the asymmetric nature of the comodule structure.
Should we have started the chain
building with $A_{01}T_0B_{01}$, the shifts would occur on the even spaces only. The transfer matrix reads:

\eb
&& t(u)=tr_0 \ \chi^t _0 A_{02N}(u,u_{2N})C_{02N-1}(u,u_{2N-1}) A_{02N-2}(u,u_{2N-2};h_{2N-1})
\ldots \non \\
&& A_{02}(u,u_2;h_{<}^{odd})C_{01}(u,u_1;h_{<}^{odd}) \ T_0(h_{<}^{odd})\
 D_{01}(u,u_1;h_{<}^{odd})
B_{02}(u,u_2;h_{<}^{odd}) \ldots \non \\
&& D_{02N-1}(u,u_{2N-1}) B_{02N}(u,u_{2N}) e^{\g \mc{D}_0}
\ee
This transfer matrix is related to the conjugate transfer matrix by:
\eb
\tilde{t}(u)=\prod_{i \textit{ even}} \chi_i^t(h_<^{odd}) t(u) \prod_{i \textit{ even}}(\chi_i^t)^{-1}(h_<^{odd})
\ee
As usual we now choose $\chi=1$.
Let us write down the transfer matrix at the point $u=u_1=0$.
\eb
&&t(0)= C_{2N,2N-1}C_{2N-2,2N-3}(h_{2N-1})\ldots C_{21}(h_{<}^{odd}) \non \\
&&P_{24}P_{46}\ldots P_{2N-2,2N}P_{2N-1,2N-3}
\ldots P_{53}P_{31} \non \\
&&B_{32}(h_1+h_{<}^{odd}) \ldots B_{2N-1,2N-2}(h_1) \non \\
&&T_{2N}(0;h_1+\ldots+h_{2N-1})P_{12N}tr_0 P_{02N}B_{02N}
e^{\g \mc{D}_0} \non
\ee
Now let us detail the last factor of this expression, i.e. the operator under the trace.
It can be written after taking into account the partial zero weight condition on $B$ in the following form.
\eb
&&tr_0( P_{02N}B_{02N}e^{\g \mc{D}_0})_{ij}=P_{klim}B_{lkmj}e^{\g \partial_k}=B_{iiij}e^{\g \partial_i}= \non \\
&&e^{\g \partial_i}B_{iiij}(\{\lambda_k +\g \delta_{ik}\})=  (e^{\g \mc{D}_{2N}})_{ii}(B_{02N}^{SC_0})_{iiij}=(e^{\g \mc{D}_{2N}}X_{2N})_{ij}
\non
\ee
where the last equality defines $X$ supposed to be invertible (such  is the case for the example given by \cite{ArCheFr}).
This condition is motivated by the fact that if we wish to
pursue the calculation of the boundary terms at some point we need its inverse.

Notice that in the bulk part of the Hamiltonian, the exponential does not appear. This is due to the preceding rewriting
of the expression $tr_0( P_{02N}B_{02N}e^{\g \mc{D}_0})$ as a product: it cancels with its inverse in $t^{-1}$ provided
the derivative does not hit $A_{02N}$ or $B_{02N}$. The final result reads.
\eb\label{HamSemi}
&&\frac{d}{du}\bigg|_{\substack{u=0\\ u_i=0 }}\log t(u)= \sum_{j=1}^N \mc{C}_{2j,2j-1}(h_{<}^{odd})+  \sum_{j=1}^{N-1} Ad\left(C_{2j+2,2j+1}(h^{odd}_{<})
\right)\cdot \mc{\check{A}}_{2j,2j+2}(h_{2j+1}+h^{odd}_{<}) +\non \\
&& \sum_{j=1}^{N-1} Ad\left(C_{2j+2,2j+1}(h^{odd}_{<}) C_{2j,2j-1}(h^{odd}_{<})\right) \cdot
\mc{B}_{2j-1,2j+2}(h_{2j+1}+h^{odd}_{<}) +\non \\
&& \sum_{j=1}^{N-1} Ad\left(C_{2j+2,2j+1}(h^{odd}_{<}) C_{2j,2j-1}(h^{odd}_{<}) B_{2j-1,2j+2}(h_{2j+1}+h^{odd}_{<})\right) \cdot
\mc{\check{D}}_{2j+1,2j-1}(h^{odd}_{<})+ \non \\
&&tr_0 \left(\check{A}'_{2N0}C_{2N,2N-1}P_{0,2N-1} B_{0,2N-1} e^{\g \mc{D}_0}\right)
X^{-1}_{2N-1} e^{-\g\mc{D}_{2N-1}} C^{-1}_{2N,2N-1}+ \non \\
&& C_{21}(h^{odd}_{<}) T_2(h_1+h^{odd}_{<}) \check{D}'_{12}(h^{odd}_{<}) T^{-1}_2(h_1+h^{odd}_{<})C^{-1}_{21}(h^{odd}_{<})+ \non \\
&& C_{2N,2N-1} e^{\g \mc{D}_{2N-1}}  X'_{2N-1} X_{2N-1}^{-1} e^{-\g \mc{D}_{2N-1}} C^{-1}_{2N,2N-1}+ \non \\
&& C_{21}(h_<^{odd})T'_2(h_1+\ldots +h_{2N-1})T_{2}(h_1+\ldots +h_{2N-1}) C_{21}^{-1}(h_<^{odd})
\ee

In the example given in \cite{ArCheFr}, the $R$-matrices contain shift operators acting on the spectral parameter. Without the inclusion
of these shift operators the matrices would obey  non standard YB equations with explicit shifts in the spectral parameters.
This will require particular caution when applying (\ref{HamSemi}) to this example.

\subsection{Fully dynamical spin chains}

The \emph{soliton preserving} conditions are again possible since the structure matrices $B$ and $C$ are of total weight zero.

For the transfer matrix
a different notation is introduced, since the two representations of the comodule both contain a shift: not only odd spaces will be shifted.
$A_{i_1,\ldots,i_m} (h_{<})$ means $A_{i_1,\ldots,i_m}$ shifted on all spaces with indices greater than the indices of $A$.

Let us write down the transfer matrix using a particular diagonal solution $\chi$ of the dual exchange relation.
\eb\label{TM22}
&&t(u)=tr_0 e^{-\g \mc{D}_0} A_{0N}A_{0N-1}(h_{<})\ldots A_{01}(h_{<}) T_{0}(h_{<}) \non \\
&&B_{01}(h_{<})\ldots B_{0N-1}(h_{<})
B_{0N} e^{\g \mc{D}_0} \chi_0^{SC_0t_0}
\ee

This transfer matrix is related to the conjugated transfer matrix as follows:
\eb
\tilde{t}(u)=\prod_i \chi_i(h_<) t(u) \prod_i \chi_i^{-1}(h_<)
\ee
For the calculation of the Hamiltonian we now suppose this redefinition has been done, and altogether eliminate $\chi_0$.

In the limit where $u \longrightarrow 0$ the permutations cancel yielding:
\eb
t(0)=n \cdot T_1(h_{<})
\ee
assumed to be invertible.
The Hamiltonian takes a very simple form:
\eb
\frac{d}{du}\bigg|_{\substack{u=0\\ u_i=0 }}\log t(u)&=& \sum_{j=1}^{N-1} \check{A}'_{j,j+1}(h_{<}) + \sum_{j=2}^{N-1}\check{B}'_{j+1,j}
(h_{<}) +T'_1(h_{<}) T^{-1}(h_{<})\non \\
&& + tr_0 e^{\g \mc{D}_0} \check{A}'_{N0
}e^{-\g \mc{D}_0}+
tr_0  e^{\g \mc{D}_0} \check{B}'_{0N}e^{-\g \mc{D}_0} +\non \\
&&  T_1(h_<) \check{B}'_{21}(h_<) T_1^{-1}(h_<)
\ee

Note that the $\exp{\g \D_0}$ seems only to contribute a shift under the trace operation. In no case can $t(u)$ in (\ref{TM22}) yield
explicit $\exp{\g \D_0}$ in its evaluation since $\chi_0$ is diagonal.

The \emph{soliton non preserving} case starts with the following transfer matrix
\eb
&&t(u)= tr_0 e^{-\g \mc{D}_0} A_{02N}(u,u_{2N})C_{02N-1}(u,u_{2N-1};h_{2N}) A_{02N-2}(u,u_{2N-2};h_{<}) \ldots \non \\
&& A_{02}(u,u_2;h_{<})C_{01}(u,u_1;h_{<}) \ T_0(h_{<}) \  D_{01}(u,u_1;h_{<})
B_{02}(u,u_2;h_{<}) \ldots \non \\
&& D_{02N-1}(u,u_{2N-1};h_{2N}) B_{02N}(u,u_{2N}) e^{\g \mc{D}_0} \chi_0^{SC_0t_0}
\ee
This transfer matrix is related to the conjugated one by the following conjugation:
\eb
\tilde{t}(u)=\prod_{i \textit{ even}} \chi_i(h_<) \ t(u) \ \prod_{i \textit{ even}} \chi_i^{-1}(h_<)
\ee
For the calculation of the Hamiltonian we use the redefined transfer matrix.
At $u=0$:
\eb
&&t(0)=C_{2N,2N-1}C_{2N-2,2N-3}(h_{<}) \ldots C_{21}(h_{<}) \non \\
&& P_{24}\ldots P_{31}P_{12N} B_{32}(h_1+h_{<})\ldots B_{2N-3,2N-2}(h_1+h_{2N-1}) B_{2N-1,2N-2}(h_1) \non \\
&& T_{1}(h_<) e^{\g \mc{D}_{2N}} X_{2N} e^{-\g \mc{D}_{2N}}
\ee
where $X_{2N}$ is defined by
\eb
X_{2N}=tr_0 P_{02N}B_{02N}
\ee
and is diagonal because of the weight zero property of $B$ as it can be shown easily:
\eb
X_{ij}=\sum_{kln}P_{klin}B_{lknj}=\sum_{kln} \delta_{kn} \delta_{il} B_{lknj} =\sum_{k} B_{ikkj}
\ee
The zero weight condition on $B$ implies that the only non-zero elements of $B$ are: $B_{iiii}, B_{ijji} (i\neq j), B_{iijj} (i\neq j)$.
Only the first and second one will give a non-zero contribution to the sum, so we have:
\eb
X_{ij}=\delta_{ij}\sum_k B_{ikki}
\ee

The Hamiltonian resulting from this transfer matrix is given by:
\eb
&&\frac{d}{du}\bigg|_{\substack{u=0\\ u_i=0 }}\log t(u)= \sum_{j=1}^N \mc{C}_{2j,2j-1}(h_{<}) +  \sum_{j=1}^{N-1} Ad\left(C_{2j+2,2j+1}(h_{<})\right)
\cdot \check{\mc{A}}_{2j+2,2j}(h_{2j+1}+h_{<})\non \\
&& \sum_{j=1}^{N-1} Ad\left(C_{2j+2,2j+1}(h_{<}) C_{2j,2j-1}(h_{<})\right) \cdot \mc{B}_{2j-1,2j+2}(h_{<}) +\non \\
&& \sum_{j=1}^{N-1} Ad\left(C_{2j+2,2j+1}(h_{<}) C_{2j,2j-1}(h_{<})
B_{2j-1,2j+2}(h_{2j+1}+h_{<}) \right) \cdot
 \mc{\check{D}}_{2j+1,2j-1}(h_{<}-h_{2j+2})+ \non \\
&& tr_0( e^{-\g \mc{D}_0}\check{A}'_{2N,0}  e^{\g \mc{D}_0} C_{2N,2N-1} P_{0,2N-1}  e^{\g \mc{D}_0} \bar{B}_{0,2N-1}) X^{-1}_{2N-1}
e^{-\g \mc{D}_{2N-1}} C^{-1}_{2N,2N-1}+\non \\
&& Ad\left(C_{21}(h_{<}) T_2(h_1+h_{<})\right) \cdot \mc{\check{D}}_{12}(h_{<}) +\non \\
&& C_{2N,2N-1} e^{\g \mc{D}_{2N-1}} X'_{2N-1}X^{-1}_{2N-1} e^{-\g \mc{D}_{2N-1}} C^{-1}_{2N,2N-1}+\non \\
&& C_{21}(h_<)T'_2(h_1+h_<) T^{-1}_2(h_1+h_<)C_{21}^{-1}(h_<)
\ee

The fifth term in this expression can be further manipulated to give an expression which clearly does not contain explicit
exponential operators.
\eb
&&tr_0( e^{-\g \mc{D}_0}\check{A}'_{2N,0}  e^{\g \mc{D}_0} C_{2N,2N-1} P_{0,2N-1}  e^{\g \mc{D}_0} \bar{B}_{0,2N-1}) X^{-1}_{2N-1}
e^{-\g \mc{D}_{2N-1}} C^{-1}_{2N,2N-1}= \non\\
&& e^{\g \mc{D}_{2N}} tr_0( \bar{\check{A}}'_{2N,0}  e^{-\g \mc{D}_{2N}} C_{2N,2N-1} e^{\g \mc{D}_{2N-1}} P_{0,2N-1}  \bar{B}_{0,2N-1})
X^{-1}_{2N-1}e^{-\g \mc{D}_{2N-1}} C^{-1}_{2N,2N-1}=\non\\
&&e^{\g \mc{D}_{2N}} tr_0(  P_{0,2N-1}\bar{\check{A}}'_{2N,2N-1}  e^{-\g \mc{D}_{2N}} C_{2N,0} e^{\g \mc{D}_{0}}  \bar{B}_{0,2N-1})
X^{-1}_{2N-1}e^{-\g \mc{D}_{2N-1}} C^{-1}_{2N,2N-1}=\non\\
&&e^{\g \mc{D}_{2N}} tr_0(  P_{0,2N-1}\bar{\check{A}}'_{2N,2N-1}  e^{\g \mc{D}_{0}} \bar{C}_{2N,0} e^{-\g \mc{D}_{2N}}  \bar{B}_{0,2N-1})
X^{-1}_{2N-1}e^{-\g \mc{D}_{2N-1}} C^{-1}_{2N,2N-1}=\non\\
&&e^{\g \mc{D}_{2N}+\g \D_{2N-1}} tr_0(  P_{0,2N-1}\bar{\check{A}}'_{2N,2N-1}(-h_0) \bar{C}_{2N,0}  \bar{B}_{0,2N-1}(-h_{2N}))
e^{-\g\D_{2N}}
e^{\g\D_{2N-1}-\g\D_{2N-1}} \non \\
&&X^{-1}_{2N-1}e^{-\g \mc{D}_{2N-1}} C^{-1}_{2N,2N-1}\non
\ee
where the notation $\bar{A}_{12}:=A^{-SL_{12}}$ was used (cf. \cite{nous}).
Note that the operator under the trace is of (triple) zero weight, thanks to the zero weight of the structure matrices. Taking the
trace on one of the spaces does not change this property on the remaining spaces, as it can be seen for example from
an argument on the equality (as sets) of incoming and outgoing indices.

This implies that the operator
\eb \non
&e^{\g \mc{D}_{2N}+\g \D_{2N-1}} tr_0(  P_{0,2N-1}\bar{\check{A}}'_{2N,2N-1}(-h_0) \bar{C}_{2N,0}  \bar{B}_{0,2N-1}(-h_{2N}))
e^{-\g\D_{2N}-\g\D_{2N-1}}&
\ee
is without explicit exponentials and so is
\eb
&e^{\g\D_{2N-1}}X^{-1}_{2N-1}e^{-\g \mc{D}_{2N-1}}&\non
\ee

Finally, we give here the results for a \emph{nondiagonal} solution $\chi$ of the dual reflection equation.
In this case the dual reflection matrix cannot be reabsorbed into the structure matrices. The main difference lies in the
nondiagonality of the object $X$ which will in turn imply the appearance of explicit exponentials in certain
boundary terms of the Hamiltonians.

In the \textit{soliton preserving} case the transfer matrix has the form:

\eb
&&t(u)=tr_0 e^{-\g \mc{D}_0} A_{0N}A_{0N-1}(h_{<})\ldots A_{01}(h_{<}) T_{0}(h_{<}) \non \\
&&B_{01}(h_{<})\ldots B_{0N-1}(h_{<})
B_{0N} e^{\g \mc{D}_0} \chi_0^{SC_0t_0}
\ee
Which gives at spectral parameter $0$:
\eb
t(u=0, u_i=0)= T_1(0, h_{<}) tr \chi^{SC}
\ee
The resulting Hamiltonian is:
\eb
\frac{d}{du}\bigg|_{\substack{u=0\\ u_i=0 }}\log t(u)&=& \sum_{j=1}^{N-1} \check{A}'_{j,j+1}(h_{<}) + \sum_{j=2}^{N-1}\check{B}'_{j+1,j}
(h_{<}) +T'_1(h_{<}) T^{-1}(h_{<})\non \\
&& tr \chi'^{SC}(tr \chi^{SC})^{-1}+ tr_0 (e^{-\g \D_0} P_{0N} A'_{N0} e^{\g \D_0} \chi_0^{SC t}) (tr \chi^{SC})^{-1} \non\\
&& T_1(h_{<}) tr_0( e^{-\g \D_0} P_{0N} B'_{0N} e^{\g \D_0} \chi_0^{SC t}) T_1^{-1}(h_{<}) (tr \chi^{SC})^{-1}+\non \\
&& T_1(h_<) \check{B}'_{21}(h_<) T_1^{-1}(h_<)
\ee

In the \textit{soliton-non-preserving} case the transfer matrix is written as:
\eb
&&t(u)= tr_0 e^{-\g \mc{D}_0} A_{02N}(u,u_{2N})C_{02N-1}(u,u_{2N-1};h_{2N}) A_{02N-2}(u,u_{2N-2};h_{<}) \ldots \non \\
&& A_{02}(u,u_2;h_{<})C_{01}(u,u_1;h_{<}) \ T_0(h_{<}) \  D_{01}(u,u_1;h_{<})
B_{02}(u,u_2;h_{<}) \ldots \non \\
&& D_{02N-1}(u,u_{2N-1};h_{2N}) B_{02N}(u,u_{2N}) e^{\g \mc{D}_0} \chi_0^{SC_0t_0}
\ee
The transfer matrix at spectral parameters $=$0 is.
\eb
&&t(0)=C_{2N,2N-1}C_{2N-2,2N-3}(h_{<}) \ldots C_{21}(h_{<}) \non \\
&& P_{24}\ldots P_{31}P_{12N} B_{32}(h_{<})\ldots B_{2N-3,2N-2}(h_{<}) B_{2N-1,2N-2}(h_<) \non \\
&& T_{1}(h_<) tr_0( P_{02N} e^{-\g \mc{D}_{2N}} B_{02N} e^{\g \mc{D}_{0}}\chi_0^{SC t})
\ee

The last factor under the trace can be rewritten as:
\eb
&& tr_0( P_{02N} e^{-\g \mc{D}_{2N}} B_{02N} e^{\g \mc{D}_{0}}\chi_0^{SC t})= e^{\g \D_{2N} } tr_0\left(P_{02N} B_{02N}^{-SL_{02N}}
\chi_0^{SCt}(-h_{2N})\right) e^{-\g\D_{2N}}= \non \\
&&e^{\g \D_{2N} }X_{2N} e^{-\g\D_{2N}} \non
\ee
where the last equality defines $X$. To go on we must suppose that $X$ is invertible. It is clear that the bulk terms of
the resulting Hamiltonian are not modified. We only give here the boundary terms.
\eb
&&H_{b}= tr_0\left(e^{-\g \D_0}\check{A}'_{2N0}e^{\g \D_0} C_{2N,2N-1} e^{\g\ D_{2N-1}} P_{02N-1}\bar{B}_{02N-1}
e^{-\g\D_{2N-1}} \chi_0^{SCt}
\right) \non \\
&&\times e^{\g\D_{2N-1}} X_{2N-1}^{-1} e^{-\g\D_{2N-1}} C_{2N,2N-1}^{-1}+ \non\\
&&C_{2N,2N-1} e^{\g \D_{2N-1}} X_{2N-1}'X_{2N-1}^{-1} e^{-\g \D_{2N-1}} C_{2N,2N-1}^{-1} + \non \\
&& C_{21}(h_<)T'_2(h_1+h_<) T^{-1}_2(h_1+h_<) C_{21}^{-1}(h_<)+\non \\
&& C_{21}(h_<)T_2(h_1+h_<)\check{D}'_{21}(h_<) T^{-1}_2(h_1+h_<) C_{21}^{-1}(h_<)
\ee

Once again note that this construction differs from that given in \cite{GZhZh}.

Those results were obtained for the $\mf{sl}_2$ case and the quasiclassical limit was
taken at different values for each of the quantum spectral parameters,
which were thus identified with non-dynamical position variables for
the spin sites. This is not the case here: all quantum spectral
parameters go to zero; the ``positions'' of the spins are supported by
the $N$ quantum space labels but are not specified and do not
enter in the dynamics at all; the $n$ dynamical variables $x_i$
are a priori unrelated to spin positions (contrary to the Ruijsenaars-
Schneider case). Their exact physical meaning is at this time an open question.

\section{Conclusion}

As a conclusion let us compute the spin chain Hamiltonian resulting from a simple example of a fully dynamical quadratic algebra.
We take the $R$-matrix which is the $\mf{gl}_2$ solution of the following dynamical Yang-Baxter equation:
\eb
&&R_{12}(\lambda+\g h_3,u_{12}) R_{13}(\lambda,u_{13}) R_{23}(\lambda+\g h_1,u_{23}) =\\
&&R_{23}(\lambda,u_{23}) R_{13}(\lambda+\g h_2,u_{13}) R_{12}(\lambda,u_{12})
\ee
where $u_{ij}=u_i-u_j$ and
\eb
&&R(\lambda,u)=E_{11}\otimes E_{11}+E_{22}\otimes E_{22}+ \alpha(\lambda, u) E_{11}\otimes E_{22} + \delta(\lambda,u) E_{22}\otimes E_{11} \non\\
&& \beta(\lambda, u) E_{12} \otimes E_{21} + \wh{\gamma}(\lambda,u) E_{21} \otimes E_{12}
\ee
with
\eb
&&\alpha(\lambda,u)=\delta(-\lambda,u)=\frac{\sinh(\lambda_{12}-\g)\sinh u}{\sinh(u-\g)\sinh \lambda_{12} }  \\
&&\beta(\lambda,u)=\wh{\g}(-\lambda,u)=\frac{\sinh(u-\lambda_{12})\sinh\g}{\sinh(u-\g)\sinh \lambda_{12}} \ee

where $\lambda_{12}=\lambda_1-\lambda_2$.

 To this Yang-Baxter equation is associated the following dynamical quadratic algebra: \eb
R_{12}(\lambda,u_1-u_2)T_1(\lambda-\g h_2,u_1) R_{21}(\lambda,u_1+u_2) T_{2}(\lambda-\g h_1,u_2) =\\
T_{2}(\lambda-\g h_1,u_2) R_{12}(\lambda,u_1+u_2) T_{1}(\lambda-\g h_2,u_1) R_{21}(\lambda,u_1-u_2)
\ee

and the corresponding dual exchange relation
\eb
R_{12}(\lambda,u_1-u_2)K_1(\lambda-\g h_2,u_1) \tilde{R}_{21}(\lambda,u_1+u_2) K_{2}(\lambda-\g h_1,u_2) =\\
K_{2}(\lambda-\g h_1,u_2) \tilde{R}_{12}(\lambda,u_1+u_2) K_{1}(\lambda-\g h_2,u_1) R_{21}(\lambda,u_1-u_2)
\ee
where
\eb
&&\tilde{R}(\lambda,u)=\tilde{\zeta}(\lambda,u)E_{11}\otimes E_{11}+\tilde{\eta}(\lambda,u)E_{22}\otimes E_{22}+ \tilde{\alpha}(\lambda, u) E_{11}\otimes E_{22} + \tilde{\delta}(\lambda,u) E_{22}\otimes E_{11} \non\\
&& \tilde{\beta}(\lambda, u) E_{12} \otimes E_{21} + \tilde{\gamma}(\lambda,u) E_{21} \otimes E_{12} \ee with \eb
&&\tilde{\zeta}(\lambda,u)=\tilde{\eta}(-\lambda,u)=\left(1-\frac{\sinh^2\g\sinh(\lambda_{12}-u)\sinh(2\g-\lambda_{12}-u)}{\sinh\lambda_{12}
\sinh(2\g-\lambda_{12})
\sinh^2(2\g-u)}\right)^{-1}\non \\
&& \tilde{\alpha}(\lambda, u)= \tilde{\delta}(-\lambda,u)=\frac{\sinh\lambda_{12} \sinh(u-\g)}{\sinh u\sinh(\lambda_{12}-\g)}\non \\
&& \tilde{\beta}(\lambda, u)=\tilde{\gamma}(-\lambda,u)=-\frac{\sinh \g \sinh\lambda_{12} \sinh(\g - u) \sinh(2 \g - \lambda_{12}
- u)} {\sinh^2(\g - \lambda_{12})\sinh(2\g - u) \sinh u} \ee

For our purpose we take the following diagonal scalar solutions: $T=\mathbbm{1} \otimes \mathbbm{1}$ and for
the dual algebra $\chi= \chi_i \, E_{ii}\otimes \mathbbm{1}$, with \eb
\chi_1(\lambda,u)=\frac{\sinh\lambda_{12}\sinh(-\lambda_1+\xi -u+\g)}{\sinh(\lambda_{12}-\g)\sinh(-\lambda_1+\xi +u-\g)} \non\\
\chi_2(\lambda,u)=\frac{\sinh\lambda_{12}\sinh(-\lambda_2+\xi-u+\g)}{\sinh(\lambda_{12}+\g)\sinh(-\lambda_2+\xi+u-\g)} \ee where $\xi$
is a free parameter\cite{YaSa}. These elements build up the following spin chain Hamiltonian: \eb
&&H=\sum_{j=1}^{N-1} h_{j,j+1}(\lambda-\g h_{<})+\non \\
&&(tr\chi^{SC})^{-1} \left\{
 f(\lambda) \bb{1} + g(\lambda) \sigma^{z}_{N}\right\} \ee

 where

 \eb h(\lambda)&=&\coth\g\left[\frac{1}{2}
\bb{1}\otimes\bb{1}-\frac{1}{2}\sigma^z \otimes \sigma^z-\sigma^{-}\otimes \sigma^{+} -\sigma^{+}
\otimes \sigma^{-}\right]+ \non \\
&&\coth\lambda_{12}\left[\frac{1}{2}\sigma^{z}\otimes \bb{1}- \frac{1}{2} \bb{1}\otimes \sigma^{z}
+\sigma^{-}\otimes \sigma^{+} -\sigma^{+}
\otimes \sigma^{-}\right] \non\\
\ee

\eb f(\lambda)&=&\frac{1}{\sinh \g}\left[\frac{\sinh(2\g -2\xi + 2\lambda_1)}{\sinh^2 (2\g-\xi+\lambda_1)}+\non
		\frac{\sinh(2\g-2\xi+2\lambda_2)}{\sinh^2(2\g -\xi +\lambda_2)}\right]+\\
&&		\frac{2}{\sinh\lambda_{12}}\left[\frac{\sinh(\g+\lambda_{12})\sinh(2\g-2\xi+2\lambda_1)}{\sinh^2(2\g-\xi+\lambda_1)}
		-\frac{\sinh(\g-\lambda_{12})\sinh(2\g-2\xi+2\lambda_2)}{\sinh^2(2\g-\xi+\lambda_2)}\right]
\ee
\eb
 g(\lambda)&=& \frac{1}{\sinh \g}\left[-\frac{\sinh(2\g -2\xi + 2\lambda_1)}{\sinh^2 (2\g-\xi+\lambda_1)}+\non
		\frac{\sinh(2\g-2\xi+2\lambda_2)}{\sinh^2(2\g -\xi +\lambda_2)}\right]+\\
	&&	\frac{1}{\sinh\lambda_{12}}\left[\frac{\sinh(\g+\lambda_{12})\sinh(2\g-2\xi+2\lambda_1)}{\sinh^2(2\g-\xi+\lambda_1)}
		+\frac{\sinh(\g-\lambda_{12})\sinh(2\g-2\xi+2\lambda_2)}{\sinh^2(2\g-\xi+\lambda_2)}\right]\\
 (tr\chi^{SC})^{-1}&=&\frac{\sinh(2\g-\xi+\lambda_1)\sinh(2\g-\xi+\lambda_2)}{2\cosh \g \sinh(\g-\xi+\lambda_2)\sinh(\g-\xi+\lambda_1)}
\ee

\section*{Acknowledgments}

The authors would like to thank Genevi\`eve Rollet and Anastasia Doikou for very helpful discussions.


\begin{thebibliography}{99}
\bibitem{ACDFR} D. Arnaudon, N. Cramp\'e, A. Doikou, L. Frappat, \'E. Ragoucy: Analytical Bethe Ansatz for closed and open $gl(\mc{N})$-spin
chains in any representation; \emph{math-ph/0411021}

\bibitem{Che} I. Cherednik; Theor. Math. Phys \tbf{61} (1984), 77.

\bibitem{KrZa}I. Krichever and A. Zabrodin, Spin Generalization of the Ruijsenaars-Schneider Model, Non-Abelian 2D Toda Chain and
Representations of Sklyanin Algebra, Uspekhi Mat. Nauk 50 (1995), n 6.\emph{hep-th/9505039}

\bibitem{Skl} E. K. Sklyanin; J. Phys A \tbf{21} (1998), 2375.

\bibitem{KuSa} P. P. Kulish, R. Sasaki; Prog. Theor. Phys. \tbf{89} (1993), 741; \emph{hep-th/9212007}

\bibitem{FrMa} L. Freidel, J. M. Maillet; Phys. Lett. B \tbf{262} (1991), 278.

\bibitem{KuSkl} P. P. Kulish, E. K. Sklyanin; J. Phys. A: Math. Gen.\tbf{24} L435 (1992)


\bibitem{KRS}P. P. Kulish, N. Yu. Reshetikhin and E. K. Sklyanin: Yang-Baxter equation and representation theory, Lett. Math. Phys.
\textbf{5} (1981), 393-403.

\bibitem{An} A. Doikou: Quantum spin chains with ``soliton non-preserving'' boundary conditions; J. Phys. A: Math. Gen. \textbf{33}
(2000) 8797-8807

\bibitem{MeNe} L. Mezincescu, R. I. Nepomechie, V. Rittenberg; Phys. Lett. A \tbf{147} (1990), 70.; L. Mezincescu, R. I. Nepomechie;
J. Phys. A: Math. Gen. \tbf{25} (1992), 2533.

\bibitem{AD} J. Avan, A. Doikou; J. Phys. A: Math. Gen.\tbf{37} (2004), 1603.

\bibitem{DM} J. Donin, A. I. Mudrov; Isr. J. Math. \tbf{136} (2003), 11.

\bibitem{DKM} J. Donin, P. P. Kulish, A. I. Mudrov; Lett. Math. Phys. \tbf{63} (2003), 179.

\bibitem{Ping} Ping Xu; Comm. Math. Phys \tbf{226} (2002), 475.

\bibitem{Et} P. Etingof, A. Varchenko: Solutions of the quantum dynamical Yang-Baxter equation and dynamical quantum groups, Comm.
Math. Phys \textbf{196} (1998), 591-640,
\emph{q-alg/9708015};B.Enriquez, P. Etingof:Quantization of classical dynamical $r$-matrices with nonabelian base \emph{math.QA/0311224};
P. Etingof: On the dynamical Yang-Baxter equation,  Proceedings of the ICM, Beijing 2002, vol. 2, 555--570 \emph{math.QA/0207008}


\bibitem{JiMiOk} M. Jimbo, T. Miwaand M. Okado, Nucl. Phys. \textbf{B300} (1998), 74.

\bibitem{GN} J. L. Gervais, A. Neveu; Nucl. Phys. B \tbf{238} (1984), 125.

\bibitem{Fe} G. Felder; Proceedings ICM Z\"urich (1994),  1247., \emph{hep-th/9407154}; Proceedings ICMP Paris (1994), 211.

\bibitem{ArCheFr} G. E. Arutyunov, L. O. Chekhov, S. A. Frolov; Comm. Math. Phys \tbf{192} (1998), 405.

\bibitem{ArFr} G. E. Arutyunov, S. A. Frolov; Comm. Math. Phys. \tbf{191} (1998), 15.

\bibitem{NAR} Z. Nagy, J. Avan, G. Rollet, Lett. Math. Phys. \tbf{67} (2004) 1.

\bibitem{NADR} Z. Nagy, J. Avan, A. Doikou, G. Rollet; \emph{math.QA/0403246}

\bibitem{FHS} H. Fan, B.-Y. Hou, K. J. Shi: Nucl. Phys. B \tbf{496} (1997), 551.

\bibitem{FHLS}H. Fan, B.-Y. Hou., G.-L. Li and K.-J. Shi:
Integrable $A_{n-1}^{(1)}$ IRF model with reflecting boundary conditions, Mod. Phys. Lett. A \textbf{26} (1997), 1929-1942

\bibitem{GZhZh} M. D. Gould, Y. Z. Zhang, S. Y. Zhao; \emph{nlin.SI/0110038}

\bibitem{KuM} P. P. Kulish, A. I. Mudrov; \emph{math.QA/0405556}

\bibitem{PPKu} P. P. Kulish: Yang-Baxter equation and reflection equations in integrable models \emph{hep-th/9507070}

\bibitem{AhnKoo} C. Ahn, W.M. Koo: Boundary Yang-Baxter equation in the RSOS/SOS representation \emph{hep-th/9508080}

\bibitem{BP} R. Behrend, P. Pearce and D. O'Brien, J. Stat. Phys \textbf{84} (1996),1.

\bibitem{ABB} J. Avan, O. Babelon, E. Billey; Comm. Math. Phys. \tbf{178} (1996), 281.

\bibitem{DeVWo} H.J. de Vega, F. Woynarovich: New integrable quantum chains combining different kinds of spins; J. Phys. A: Math. Gen
\textbf{25} (1992) 4499-4516.

\bibitem{nous} Z. Nagy, J. Avan, A. Doikou, G. Rollet: Commuting quantum traces for quadratic algebras; \emph{math.QA/0403246}

\bibitem{YaSa} W.-L. Yang, R. Sasaki: Solution of the dual reflection equation for $A_{n-1}^{(1)}$ solid-on-solid model;
J. Math, Phys \tbf{45} (2004), 4301.

\bibitem{BaDo} P. Baseilhac, K. Koizumi: Sine-Gordon quantum field theory on the half-line with quantum boundary degrees of freedom; 
Nucl. Phys. \tbf{B649} (2003) p. 491-510; \emph{hep-th/0208005}

\end{thebibliography}
\end{document}